\documentclass[11pt]{article}
\usepackage{bbm}
\usepackage{mathrsfs}
\usepackage{amscd}
\usepackage{amsmath,amsfonts,amssymb,amscd}
\usepackage{indentfirst,graphics,epsfig,psfrag}
\input{epsf}
\usepackage{ifpdf}
\usepackage{enumerate}
\usepackage{appendix}
\usepackage{enumerate}
\usepackage{amsmath}
\usepackage{lineno}
\usepackage{color}
\makeatletter \@addtoreset{figure}{section} \makeatother
\makeatletter
\long\def\@makecaption#1#2{%
   \vskip 10\p@
   \setbox\@tempboxa\hbox{{#1}\ \ #2}%
   \ifdim \wd\@tempboxa >\hsize

       {#1}\ \ #2\par
   \else
       \hbox to\hsize{\hfil\box\@tempboxa\hfil}%
   \fi}
\makeatother

\newtheorem{thm}{Theorem}[section]
\newtheorem{cor}[thm]{Corollary}
\newtheorem{lem}[thm]{Lemma}

\newtheorem{exe}[thm]{Example}
\newtheorem{obs}[thm]{Observation}
\newtheorem{pro}[thm]{Proposition}

\newcommand{\qed}{{\hfill\rule{3pt}{7pt}}}

\def\qed{\hfill \rule{4pt}{7pt}}

\begin{document}
\title{\textbf{On the regular $k$-independence number of graphs} \footnote{Research supported by the National Science Foundation of
China (Nos. 61164005, 11161037, 11101232, 61460005 and 11461054), the National
Basic Research Program of China (No. 2010CB334708) and the Program
for Changjiang Scholars and Innovative Research Team in Universities
(No. IRT1068), the Research Fund for the Chunhui Program of Ministry
of Education of China (No. Z2014022) and the Nature Science
Foundation from Qinghai Province (Nos. 2012-Z-943, 2014-ZJ-721 and
2014-ZJ-907).}}
\author{
\small Zhiwei Guo$^1$, \ Haixing Zhao$^2$, \ Hong-Jian Lai$^3$, \ Yaping Mao$^1$\\[0.3cm]
\small $^1$Department of Mathematics, Qinghai Normal\\
\small University, Xining, Qinghai 810008, China\\
\small E-mails: guozhiweic@yahoo.com; maoyaping@ymail.com.\\[0.3cm]
\small $^2$School of Computer, Qinghai Normal\\
\small University, Xining, Qinghai 810008, China\\
\small E-mail: h.x.zhao@163.com\\[0.3cm]
\small $^3$Department of Mathematics, West Virginia\\
\small University, Morgantown, WV 26506, USA\\
\small E-mail: hjlai2015@hotmail.com\\[0.3cm]}
\date{}
\maketitle
\begin{abstract}
The \emph{regular independence number}, introduced by Albertson and
Boutin in 1990, is the maximum cardinality of an independent set of
$G$ in which all vertices have equal degree in $G$. Recently, Caro,
Hansberg and Pepper introduced the concept of regular $k$-independence
number, which is a natural generalization of the regular
independence number. A \emph{$k$-independent set} is a set of
vertices whose induced subgraph has maximum degree at most $k$.
The \emph{regular $k$-independence number} of $G$, denoted by $\alpha_{k-reg}(G)$,
is defined as the
maximum cardinality of a $k$-independent set of $G$ in which
all vertices have equal degree in $G$. In this paper, the exact
values of the regular $k$-independence numbers of some special
graphs are obtained. We also get some lower and upper bounds
for the regular $k$-independence number of trees with given
diameter, and the lower bounds for the regular $k$-independence number of line graphs.
For a simple graph $G$ of order $n$, we show that
$1\leq\alpha_{k-reg}(G)\leq n$ and characterize the extremal graphs.
The Nordhaus-Gaddum-type results for the regular $k$-independence
number of graphs are also obtained.
\\[2mm]

\\[2mm]

\end{abstract}

\section{Introduction}

Graphs considered in this paper are undirected, finite and simple. We refer
to \cite{bondy} for undefined notations and terminology.
In particular,
we use
$L(G)$, $\bar{G}$, $\Delta(G)$ and $\delta(G)$ to denote the
line graph, the complementary graph, the maximum degree and
minimum degree of a graph $G$,
respectively. If $X \subseteq V(G)$ or $X \subseteq E(G)$, then
$G[X]$ is the subgraph of $G$ induced by $X$. For integers $i,j \ge 0$, let $D_i(G)$
denote the set of degree $i$ vertices of $G$, and $D_{\ge i}(G) = \bigcup_{j \ge i}D_j(G)$.
A subset $X \subseteq V(G)$ is \emph{regular}
if for some $i$ with $\delta(G) \le i \le \Delta(G)$, $X \subseteq D_i(G)$;
and is independent if $\Delta(G[X]) = 0$.
An independent set $X$ of $G$ is  a \emph{regular independent set} if
$X$ is also regular.

The \emph{regular independence number}, denoted $\alpha_{reg}(G)$
and introduced by Albertson and Boutin {\upshape \cite{Albertson}}
in 1990, is defined to be the maximum cardinality of an independent
set of $G$ in which all vertices have equal degree in $G$. The
parameter $\alpha_{reg}(G)$ is closely related to fair domination
number $fd(G)$ introduced in  {\upshape \cite{Henning}}. A
\emph{fair dominating set} is a set $S\subseteq V(G)$ such that all
vertices $v \in V (G)\setminus S$ have exactly the same non-zero
number of neighbors in $S$. The \emph{fair domination number}
$fd(G)$ is the cardinality of a minimum fair dominating set of $G$.
By definition, if $\delta(G)\ge 1$ and $R$ is a maximum regular
independent set of $G$, then $V(G)\setminus R$ is a fair dominating
set of $G$. A vertex subset $S \subseteq V(G)$ of $G$ is
\emph{$k$-independent set} if $\Delta(G[S]) \le k$. The
\emph{$k$-independence number}, denoted $\alpha_{k}(G)$, as the
maximum cardinality of a $k$-independent set. There have been quite
a few studies on $k$-independent sets, as seen in
\cite{Tuza,Hansberg,Pepper}, among others. For
$k$-independent set and $k$-independent number, Chellali, Favaron,
Hansberg, and Volkmann published a survey paper on this subject; see
\cite{Chellali}.

Recently, Caro, Hansberg and Pepper {\upshape \cite{Han}}
introduced the concept of regular $k$-independence number, which
naturally generalizes both
the regular independence number and the $k$-independence number.
The \emph{regular $k$-independence number} of a graph $G$, denoted $\alpha_{k-reg}(G)$,
is defined to be
the maximum cardinality of a regular $k$-independent set of $G$.
More precisely, for nonnegative integers $k$ and $j$,
we define $\alpha_{k,j}(G) = \max\{|X|: $
$X$ is a $k$-independent set of $G$ and $X \subseteq D_j(G) \}$.
It follows by definition that,

\begin{equation} \label{a-kj}
\alpha_{k-reg}(G)=\max \{\alpha_{k,j}(G), \delta(G) \leq j \leq\Delta(G) \}.
\end{equation}
When $k=0$, $\alpha_{0-reg}(G)=\alpha_{reg}(G)$ and for regular graphs, $\alpha_{reg}(G)=\alpha(G)$ and $\alpha_{k-reg}(G)=\alpha_{k}(G)$.

For each integer $i \ge 0$, define $n_i(G) = |D_i(G)|$. We often write
$n_i$ for $n_i(G)$ when the graph $G$ is understood from the context.
The \emph{repetition number} of $G$, denoted $rep(G)$, was introduced in
{\upshape \cite{Wes}} and defined as the maximum number of vertices
with equal degree in $G$. Thus
\begin{equation} \label{rep}
rep(G) = \max\{|D_i(G)|: \delta(G) \le i \le \Delta(G)\}.
\end{equation}

The notation of $\chi_k(G)$ is the \emph{$k$-chromatic number}
of $G$, defined as the minimum number of colors needed to color
the vertices of the graphs $G$ such that the graphs induced by the
vertices of each color class have maximum degree at most $k$.
Note that $\chi_0(G)$ is the
classic
chromatic number $\chi(G)$.

In {\upshape \cite{Han}}, Caro, Hansberg and Pepper investigated the
regular $k$-independence number of trees and forests, and they
generalized and extended the results of Albertson and Boutin
{\upshape \cite{Albertson}} in to $\alpha_{k-reg}(G)$. They
presented a lower bound on $\alpha_{k-reg}(G)$ for $k$-trees and
gave analogous results for $k$-degenerate graphs and some specific
results about planar graphs, and then gave lower bounds on
$\alpha_{2-reg}(G)$ for planar and outerplanar graphs. The authors
also analyzed complexity issues of regular $k$-independence.

This paper is organized as follows. In Section $2$,
the exact values of the regular $k$-independence numbers of
complete graphs, complete multipartite graphs, paths, cycles
and stars are
determined.
Sharp bounds for the regular $k$-independence number
of trees with given diameter
are obtained in Section $3$.
In Section $4$, we obtain the lower bounds for the regular
$k$-independence number of general $m$-vertex line graphs.
For some families of sparse graphs  such as trees,
maximal outerplanar graphs and triangulations,
we present lower bounds for the regular $k$-independence
number of their line graphs, which improve several former results.
For a simple graph $G$ of order $n$, we show that
$1\leq\alpha_{k-reg}(G) \leq n$,
and characterize all extremal graphs in Section $5$.

Let $\mathcal {G}(n)$ denote the class of simple graphs of order $n
\ (n\geq 2)$. For $G\in \mathcal {G}(n)$, $\Bar{G}$ denotes the complement
of $G$. Give a graph parameter $f(G)$ and a positive integer $n$, the
\emph{Nordhaus-Gaddum(\textbf{N-G}) Problem} is to determine sharp
bounds
 for both $f(G)+f(\Bar{G})$ and  $f(G)\cdot
f(\Bar{G})$,
as $G$ ranges over the class $\mathcal {G}(n)$,
and characterize the extremal graphs.
The Nordhaus-Gaddum type relations have received wide
attention,
as seen in the survey paper \cite{Aouchiche} by Aouchiche
and Hansen. The Nordhaus-Gaddum-type problem
on the regular $k$-independence number of graphs
is studied in Section $6$.

\section{Results for some special graphs}

In this section, we will determine the regular $k$-independence numbers
in several special families of graphs. Throughout this section,
$n > 0$ denotes an integer.

\begin{pro}\label{pro2-1}
Let $K_n$ be a complete graph of order $n$. Then
$$
\alpha_{k-reg}(K_n)=
\begin{cases}
i+1,&~if~ k=i~(0\leq i\leq n-1);\\
n,  &~if~ k\geq n.
\end{cases}
$$
\end{pro}
\begin{pf}
Let $K_n$ be a complete graph with $n$ vertices. Since $K_n$ is a
regular graph, it follows that
$\alpha_{k-reg}(K_n)=\alpha_{k}(K_n)$. By the definition of
$\alpha_{k}(K_n)$, if $k = i$ for some $i$ with $0 \le i \le n-1$,
then every $i+1$ subset of vertices is a maximum regular
$k$-independent set of $K_n$; if $k \ge n$, then $V(K_n)$ is a
maximum regular $k$-independent set of $K_n$. This proves the
proposition. \qed
\end{pf}

\begin{pro}\label{pro2-2}
Let $K_{r_1,r_2,\cdots,r_n}$ be a complete $n$-partite graph.
\\
$(1)$ If $r_1=r_2= \cdots=r_n=a$, then
$$
\alpha_{k-reg}(K_{r_1,r_2,\cdots,r_n})=
\begin{cases}
ia,   &~if~ (i-1)a\leq k < ia~(1\leq i\leq n); \\
na,  &~if~ k\geq na.
\end{cases}
$$
$(2)$ If $r_1<r_2< \cdots<r_n$, then $\alpha_{k-reg}(K_{r_1,r_2,\cdots,r_n})=r_n$ for $k\geq 0$.
\\
$(3)$ If $r_1<\cdots< r_i=r_{i+1}=\cdots =r_j<r_{j+1}< \cdots <r_n~(i<j)$, then
$$
\alpha_{k-reg}(K_{r_1,r_2,\cdots,r_n})=
\begin{cases}
\max\{mr_i, ~r_n\},       &~if~ (m-1)r_i\leq k < mr_i~(1\leq m \leq (j-i));\\
\max\{(j-i)r_i, ~r_n\} ,  &~if~ k\geq (j-i)r_i.
\end{cases}
$$
\end{pro}

\begin{pf}
Let $G= K_{r_1,r_2,\cdots,r_n}$
with partite sets $V_1, V_2, \cdots, V_n$
such that
$|V_j| = r_t$, $1 \le t \le n$.

\noindent $(1)$
Assume that $r_1=r_2= \cdots=r_n=a$. Then
$K_{r_1,r_2,\cdots,r_n}$ is a regular graph, and so $\alpha_{k-reg}(K_{r_1,r_2,\cdots,r_n})=\alpha_{k}(K_{r_1,r_2,\cdots,r_n})$.
For each $i \in \{1, 2, \cdots, n\}$ with $(i-1)a\leq k < ia$,
any union of $i$ of the partite sets is a maximum regular $k$-independent set;
and for $k \ge na$, $V(G)$ is the only maximum regular $k$-independent set.
This justifies (1).

\noindent
$(2)$
Assume that $r_1<r_2< \cdots<r_n$. By the definition of regular independent sets
and since the $r_i$'s are mutually distinct,
a vertex subset $X$ of $G$ is a regular independent set
if and only if $X \subseteq V_t$ for
some $t$ with $1 \le t \le n$. It follows that
$V_n$ is the only maximum
regular $k$-independent set of $G$.
Thus $\alpha_{k-reg}(K_{r_1,r_2,\cdots,r_n})=r_n$ for $k \geq 0$.

\noindent
$(3)$
Assume that for some integers $1 \le i < j \le n$,
$r_1<\cdots< r_i=r_{i+1}=\cdots =r_j<r_{j+1}< \cdots <r_n$.
For a fixed integer $k$, let $X$ be a maximum regular $k$-independent set of $G$.
By assumption of (3) and by the definition of regular independent sets,
we note that
either $X \subseteq V_t$ for
some $t$ with $1 \le t \le n$ or $X \subseteq \bigcup_{t=i}^j V_t$.
If for an integer $m$ with $1 \leq m \leq (j-i)$,
we have $(m-1)r_i\leq k < mr_i$, then either $X \subseteq \bigcup_{t=i}^j V_t$
with $|X| = mr_i$
or $X = V_n$, whence $|X| = \max\{mr_i, ~r_n\}$.
If $k\geq (j-i)r_i$, then either $X = \bigcup_{t=i}^j V_t$
or $X = V_n$, whence $|X| = \max\{(j-i)r_i, ~r_n\}$.
This verifies (3).
\qed
\end{pf}

\begin{pro}\label{pro2-5}

Let $m\geq 2$ be an integer, $i \in \{0,1,2\}$,
$P_n$ be a path of order $n$, where $n=3(m-2)+2+i$.
Each of the following holds.

\noindent
$(1)$ If $m\geq 3$, then
$$
\alpha_{k-reg}(P_n)=
\begin{cases}
\lceil \frac{n-2}{2}\rceil, &~if~ k=0; \\
n-m,  &~if~ k=1;\\
n-2,  &~if~ k\geq 2.
\end{cases}
$$
\noindent
$(2)$ If $m=2$, then
$$
\alpha_{k-reg}(P_2)=
\begin{cases}
1, &~if~ k=0; \\
2,  &~if~ k\geq 1.
\end{cases}
$$
and
$$
\alpha_{k-reg}(P_3)=\alpha_{k-reg}(P_4)=
2
$$
for $k\geq 0$.
\end{pro}

\begin{pf}
As the proof for (2) is straightforward, we only need to
show the validity of (1).
Assume that $m\geq 3$. Then $n\geq 5$, and $D_2(P_n) = n-2$.
If $k=0$, then $D_2(P_n)$ contains an independent subset $W_0$ with
$|W_0| = \lceil \frac{n-2}{2}\rceil$.
Likewise, if $k=1$, then $D_2(P_n)$ contains a 1-independent set $W_1$
consisting $i$ isolated vertices and the vertex set of a matching with
$\frac{n-2-i}{3}$ edges. It is routine to show that $W_1$ is a maximum
1-independent set, and so
$\alpha_{1-reg}(P_n)=2(\frac{n-2-i}{3})+i=2(m-2)+i=2m-4+i=n-m$.
If $k\geq 2$, then $D_2(P_n$ is a maximum
regular $k$-independent set, and so $\alpha_{k-reg}(P_n)=n-2$.

\qed
\end{pf}

\begin{pro}\label{pro2-6}
Let $C_n$ be a cycle of order $n$. Then
$$
\alpha_{k-reg}(C_n)=
\begin{cases}
\lfloor\frac{n}{2}\rfloor,       &~if~ k=0;\\
2a, &~if~ k=1~and~n=3a~or~n=3a+1;\\
2a+1, &~if~ k=1~and~n=3a+2;\\
n,  &~if~k\geq 2.
\end{cases}
$$
\end{pro}
\begin{pf}
Denote $C_n = v_1v_2 \cdots v_nv_1$. As $C_n$ is regular, any
$k$-independent set of $C_n$ is also a regular $k$-independent set.
If $k=0$, then $\{v_{2i+1}: 0 \le i \le \frac{n}{2}-1\}$ is a
maximum regular independent set of $C_n$. Hence
$\alpha_{0-reg}(C_n)= \lfloor \frac{n}{2}\rfloor$. If $k=1$, then
$V(C_n) - \{v_{3i}: 1 \le i \le  \frac{n}{3} \}$ (if  $n \equiv 0$
(mod 3)) or $V(C_n) - (\{v_{3i}: 1 \le i \le
\lfloor\frac{n}{3}\rfloor\} \cup \{v_n\})$ (if  $n \equiv 1$ or $n
\equiv 2$ (mod 3)) is a maximum regular $k$-independent set of
$C_n$. If $k \ge 2$, $V(C_n)$ is a maximum regular $k$-independent
set. This justifies the proposition. \qed
\end{pf}

\begin{pro}\label{pro2-7}
Let $S_{1,n-1}$ be a star of order $n$. Then $\alpha_{k-reg}(S_{1,n-1})=n-1$ for $k\geq 0$.
\end{pro}
\begin{pf}
Since $S_{1,n-1}$ has $n$ vertices, there are $n-1$ vertices of degree $1$ in $S_{1,n-1}$. When $k\geq 0$, the subgraph induced by all the vertices of degree $1$ in $S_{1,n-1}$ is a regular $k$-independent set. By the definition of the regular $k$-independence number, $\alpha_{k-reg}(S_{1,n-1})=n-1$. \qed
\end{pf}

\section{Results for trees with given diameter}

Caro, Hansberg and Pepper {\upshape \cite{Han}} generalized and
extended the result that $\alpha_{reg}(T)\geq\frac{n+2}{4}$ for any
tree $T$, obtained by Albertson and Boutin in {\upshape
\cite{Albertson}. They showed that for every tree $T$ on $n\geq 2$
vertices, $\alpha_{k-reg}(T)\geq\frac{2(n+2)}{7}$ for $k=1$ and
$\alpha_{k-reg}(T)\geq\frac{(n+2)}{3}$ for $k\geq 2$. In this
section, we improve the bound of $\alpha_{k-reg}(T)$ for $k\geq 2$
by considering the diameter of a tree.} Throughout this section, let
$n \ge 8$, $t > 0$ and $k \ge 2$ be integers, and $T_{n,t}$ denote
the family of trees with order $n$ and diameter $n-t$. For
notational convenience, when it is clear form the context, we also
use $T_{n,t}$ to denote a member in this family. Thus $T_{n,1}$ is a
path with $n$ vertices and $T_{n, n-2}$ is a star with $n$ vertices.
As the regular $k$-independence number of paths and stars are
determined in the section above, we always assume that $2 \leq t\leq
n-3$ in this section. If $G = T_{n,t}$, then for each $i \ge 1$, let
$n_i = |D_i(G)|$ and $N_3 = |D_{\ge 3}(G)|$. In a graph $H$, an
\emph{elementary subdivision} of an edge $uv$ is the operation of
replacing $uv$ with a path $uwv$ through a new vertex $w$. A
\emph{subdivision} of $H$ is the graph obtained by a finite sequence
of elementary subdivisions on $H$. As usual, a \emph{leaf} of a tree
is a vertex of degree $1$ in the tree. The main purpose of this
section is to investigate the regular independence number for trees
with given diameters. We start with lemmas and examples.

\begin{lem} \label{n1}
Let $T$ be a tree on $n \ge 2$ vertices. Each of the following holds.
\\
(i) $|D_1(T)| = \sum_{v \in D_{\ge 3}(G)} (d_T(v) - 2) + 2$.
\\
(ii) If $T = T_{n,t}$, then $N_3 = |D_{\ge 3}(T_{n,t})| = n - n_1-n_2$.
and $n_1 \ge N_3+2$.
\end{lem}
\begin{pf} We outline our proofs. Lemma \ref{n1}(i) holds if
$|V(T)| = 2$, and so it can be justified  by induction on $|V(T)|$,
by considering $T - v$ for some $v \in D_1(T)$ in the inductive
step. Lemma \ref{n1}(ii) follows from the definitions.\qed
\end{pf}

\begin{exe} \label{Example1}
Let $\ell_1 \ge \ell_2 \cdots \ell_r \ge 1$ and $r \ge 3$ be integers, and let
$K_{1,r}$ denote the tree with a vertex $v_0$ of degree $r$ and
$D_1(K_{1,r})=\{v_1, v_2, \cdots, v_r\}$.
\\
(i)
Define $K_{1,r}(\ell_1, \ell_2, \cdots, \ell_r)$ to be the graph obtained from
$K_{1,r}$ by replacing each edge $v_0v_i$ by a $(v_0,v_i)$-path of order $\ell_i$,
for each $i$ with $1 \le i \le r$. When $\ell_3 = \ell_{r-1}$, we also use
$K_{1,r}(\ell_1, \ell_2, \ell_3^{r-3}, \ell_r)$ for $K_{1,r}(\ell_1, \ell_2, \cdots, \ell_r)$.
Let $T = K_{1,r}(\ell_1, \ell_2, \ell_3^{r-3}, \ell_r)$.
It is elementary to compute that $T = T_{n,t}$ with $n = 1 + \sum_{i=1}^r (\ell_i - 1)$
and $n  - t = \ell_1 + \ell_2-2$. If $k \ge 2$, then
\[
\alpha_{k-reg}(T) =
\left\{
\begin{array}{ll}
r & \mbox{ if $r > \frac{n-1}{2}$}
\\
n - r - 1 & \mbox{ if $r \le \frac{n-1}{2}$}
\end{array} \right. .
\]
(ii) Let $T$ be a tree with $|D_{\ge 3}(T)| \ge 2$. Assume that $z, z' \in D_{\ge 3}(T)$
and $T_1, T_2$ be the two subtrees of $T$ such that
$T = T_1 \cup T_2$, $V(T_1) \cap V(T_2) = \{z'\}$, $z \in V(T_1)$ and $z' \in D_2(T_1)$.
View $T_2'$ is a copy of $T_2$ but vertex disjoint from $V(T_1)$.
Obtain a new tree $T'$ from the vertex disjoint union of $T_1$ and $T_2'$ by identifying
$z \in V(T_1)$ and $z' \in V(T_2')$. We use $O_{z' \rightarrow z}$ to denote
this operation and write $T'=O_{z' \rightarrow z}(T)$, and
use $O_{z' \leftarrow z}$ to denote the reverse operation. Hence $T= O_{z' \leftarrow z}(T')$.
By definition and by Lemma \ref{n1},
$|V(T)| = |V(T')|$, $|D_1(T)|  = |D_1(T')|$ and $|D_{\ge 3}(T)| - 1 = |D_{\ge 3}(T')|$.
For a fixed $z$, define relation $T \sim T'$ if and only if for some $z'$, $T'=O_{z' \rightarrow z}(T)$. Then
$\sim$ is an equivalence relation on the set of all trees with the same number of vertices and same number of
leaves.
\\
(iii) Let $T$ be a given tree $T$ with $|D_{\ge 3}(T)| \ge 1$, and
let $z \in D_{\ge 3}(T)$ be a fixed vertex. Define
${\cal F}(T, z)$ to be equivalence class containing $T$ under the relation $\sim$ defined in
(ii) above. By definition, $T \in {\cal F}(T, z)$.
\\
(iv)
Suppose that $n$ and $t$ are integers with $2 \le t \le n-3$.
Let $h = \lfloor \frac{n-t}{2} \rfloor$. Since $t \ge 2$, we can
write $n-1 = qh + r$ for some integers $q \ge 2$ and $1 \le r \le h$.
Define
\[
T(n,t) = \left\{\begin{array}{ll}
K_{1, q+1}(h+1, h+1, (h+1)^{q-2}, r+1) & \mbox{ if $n - t \equiv 0$ (mod 2)}
\\
K_{1, q+1}(h+2, h+1, (h+1)^{q-2}, r+1) & \mbox{ if $n - t \equiv 1$ (mod 2)}
\end{array} \right. ,
\]
and
let $z_0$ be the only vertex in $T(n,t)$ with degree $q+1$.
By Example \ref{Example1}(iii),  for each $T \in {\cal F}(T(n,t),z_0)$,
$|D_1(T)| = q+1$. By definition, the diameter of $T(n,t)$ is $n-t$.
Direct computation yields that
$|D_1(T(n,t))|=q + 1$ and
$|D_2(T(n,t))| = n - q - 2$.
\end{exe}

For integers $n > t \geq 2$ with $t \le n-3$, define
\begin{equation} \label{ff}
f(n,t) =
\left\{
\begin{array}{ll}
\displaystyle \left\lceil \frac{2(t-1)}{n-t} \right\rceil + 2 &
\mbox{ if $n - t \equiv 0$ (mod 2)}
\\[10pt]
\displaystyle \left \lceil \frac{2(t-1)}{n-t-1} \right\rceil + 2 &
\mbox{ if $n - t \equiv 1$ (mod 2)}
\end{array} \right. ,
\end{equation}

\begin{lem} \label{n1-bound}
Suppose that $T = T_{n,t}$ with $2 \le t \le n-3$.
Let  $P = v_1v_2, \cdots, v_{n-t+1}$ be a longest path
in $T(n,t)$ (as defined in Example \ref{Example1}(iv)),
$h = \lfloor \frac{n-t}{2} \rfloor$,  and {\color{blue}$z_0 = v_h$} {\color{red}$z_0 = v_{h+1}$}.
Express $n-1 = qh + r$ for some integers $q \ge 2$ and $1 \le r \le h$.
Then each of the following holds.
\\
(i) $n_1 = |D_1(T)| \ge f(n,t)$.
\\
(ii) Equality in (i) holds if and only if both $q + 1 = f(n,t)$ and
$T \in {\cal F}(T(n,t),z_0)$.
\end{lem}

\begin{pf}

Since $T$ is connected and since
$N_3 > 0$, there must be a $j_0$ with $v_{j_0} \in V(P) \cap D_{\ge 3}(G)$. Without lose
of generality, we assume that $1 < j_0 < n-t+1$
such that $|\lceil \frac{n-t+1}{2} \rceil - j_0|$ is minimized.
By symmetry, we may assume that $1 < j_0 \le \lceil \frac{n-t+1}{2} \rceil = h+1$
We shall argue by induction on $N_3$. Since $t \le n-3$,  we have
$N_3 > 0$.

Suppose that $N_3 = 1$.
Then for any $w \in D_1(T) - \{v_1, v_{n-t+1}\}$,
there exists a unique $(w, v_{j_0})$-path $P_w$ in $T$ such that
$V(P_w) \cap V(P) = \{v_{j_0}\}$.

Assume first that $n - t \equiv 0$ (mod 2), and so $n - t = 2h$.
Since the diameter of $T$ is $n - t$, and since $j_0 \le h+1$ for
any $w \in D_1(T) - \{v_1, v_{n-t+1}\}$, $|E(P_w)| \le j_0- 1\le h$.
It follows that
\begin{eqnarray} \label{n-t=0}
n - 1 & = & |V(T) - \{v_{j_0}\}| = |V(P) - \{v_{j_0}\}| + \sum_{w \in D_1(T) - \{v_1, v_{n-t+1}\}} |V(P_w - v_{j_0})|
\\ \nonumber
& \le & |V(P)| - 1+ (n_1 - 2) (j_0-1)
\le n - t + (n_1-2)h,
\end{eqnarray}
and so $n_1 \ge f(n,t)$. Assume that we have $n_1 = f(n,t)$.
Then, if $h$ divides $t-1$, then  every
inequality in (\ref{n-t=0}) must be an equality; and
if $h$ does not divide $t-1$, then for some integer $r'$
with $0 < r' < h$,
$n - 1 = n - t + (n_1-2)h - r'$. It follows that $j_0 = h+1 = |V(P_w)|$,
for all but at most one $w \in D_1(T) - \{v_1, v_{n-t+1}\}$.
Since
$n - 1 = qh + r$ with $1 \le r \le h$, we have  $n_1 = q+1$.
As $N_3= 1$, $T$ must be a subdivision of $K_{1, n_1}$, and so
$T =  K_{1, n_1}(h+1, h+1, (h+1)^{q-2}, r+1) = T(n,t)$.

The proof for the case when  $n - t \equiv 1$ (mod 2)
is similar, using $n - t = 2h + 1$ and $(n_1-2)h + n-t \ge n-1$ instead, and so
it is omitted.

We now assume that $N_3 > 1$, and that Lemma \ref{n1-bound} holds
for smaller values of $N_3$. Since $N_3 \ge 2$, there exists a $w \in D_{\ge 3}(T) - \{v_{j_0}\}$.
Let $T' = O_{w \rightarrow v_{j_0}}$.
By Example \ref{Example1},
$|D_1(T')| = |D_1(T)| = n_1$. As $j_0$ is so chosen that $|\lceil \frac{n-t+1}{2} \rceil - j_0|$
is minimized, the diameter of $T'$ is also $n - t$. However,
$D_{\ge 3}(T') = D_{\ge 3}(T) - \{w\}$. By induction,
\[
n_1 = |D_1(T)| = |D_1(T')| \ge f(n,t).
\]
If equality holds, then by induction, $T' \in {\cal F}(T(n,t),v_{j_0})$, where
$j_0 = h+1$. This complete the proof of the lemma.
\qed
\end{pf}

\begin{lem} \label{max=}
Suppose that $k \ge 2$ and $T = T_{n,t}$ with $2 \le t \le n-3$. %If $\frac{n}{3} < t < n-5$,
Then $\alpha_{k-reg}(T) =\max\{|D_1(T)|, |D_2(T)|\}$.
\end{lem}

\begin{pf} Since $k \ge 2$, both $D_1(T)$ and $D_2(T)$ are regular
2-independent sets of $T$. Therefore, $\alpha_{k-reg}(T)  \ge \max\{|D_1(T)|, |D_2(T)|\}$.
If $X$ is a maximum k-independent set of $T$, then for some
$i$, $X \subseteq D_i(T)$, and so $\alpha_{k-reg}(T) = |X| \le
|D_i(T)| \le \max\{|D_i(T)|: i \ge 1\}$. By Lemma \ref{n1}(ii),
$|D_1(T)| \ge N_3 + 2 > N_3 = \sum_{j \ge 3} |D_j(T)|$. This implies that
$\max\{|D_1(T)|, |D_2(T)|\} \le \alpha_{k-reg}(T)  \le \max\{|D_i(T)|: i \ge 1\} \le \max\{|D_1(T)|, |D_2(T)|\}$.
\end{pf}\qed

\begin{thm}\label{th3-1}

Let $k \geq 2$ be an integer and $T_{n,t}$ be a tree with
order $n \ge 8$ and diameter $n-t$ with $2 \le t \le n-3$.

\noindent
(i) If $2\leq t\leq\frac{n-1}{3}$, then
\begin{equation} \label{3-1a}
n-2t\leq\alpha_{k-reg}(T_{n,t})\leq n-4.
\end{equation}
\noindent

\noindent
(ii)  If $\frac{n}{3}\leq t \leq n-5 \ (n\geq 8)$, then
\begin{equation} \label{3-1b}
\frac{n+2}{3}\leq\alpha_{k-reg}(T_{n,t})\leq
\max\left\{n-f(n,t)-1,t+1\right\}.
\end{equation}

\noindent
(iii) If $t=n-4$, then
\begin{equation} \label{3-1c}
\left\lceil\frac{n-1}{2}\right\rceil\leq\alpha_{k-reg}(T_{n,t})\leq
t+1.
\end{equation}

\noindent (iv) If $t=n-3$, then
\begin{equation} \label{3-1d}
\alpha_{k-reg}(T_{n,t})=t+1.
\end{equation}

\end{thm}

\begin{pf}
%Let $T = T_{n,t}$.
Since $k\geq 2$, both $D_1(T_{n,t})$ and $D_2(T_{n,t})$
are  regular $k$-independent set.
By the definition of $T_{n,t}$,
$n_1 \le t+1$.
By Lemma \ref{n1-bound} and by $N_3 \ge 1$,

$n_2 \le n - f(n,t) - 1$,

where equality holds if and
only if $N_3 = 1$ and $n_1 = f(n,t)$. By Lemma \ref{max=},
\begin{equation} \label{u-bound}
\alpha_{k-reg}(T_{n,t}) =\max\{n_1, n_2\}
\le \max\{ n-f(n,t)-1, t+1\}.
\end{equation}

\noindent
(i) Suppose $2\leq t\leq\frac{n-1}{3}$. Since $n_2=n-n_1-N_3$
and $n_1\geq N_3+2$, it follows that $n_2\geq n-2n_1+2$.
Since $diam(T_{n,t})=n-t$, we have $n_1 \le t+1$.
If $t\leq\frac{n-1}{3}$, then  $n_2 \geq n_1$.
Thus $D_2(G)$ is a
maximum  regular $k$-independent set
and so $\alpha_{k-reg}(T_{n,t})=n_2$.
By definition, $T_{n,t}$ contains a path $P=v_1v_2\cdots v_{n-t+1}$.
As the remaining $t-1$ vertices in $V(T_{n,t}) - V(P)$
are adjacent to at most $t-1$ vertices in $D_2(P)$,
it follows that $n_2 \ge |V(P)| - 2 - (t-1) = n-2t$.

Since $2\leq t\leq\frac{n-1}{3}$, we have $t-1\leq \lfloor\frac{n-t}{2}\rfloor$.
As $T_{n,t}$ is a tree, the remaining $t-1$ vertices in $V(T_{n,t}) - V(P)$
are adjacent to at least one vertex in $D_2(P) \cap D_{\ge 3}(G)$ and contains at least one
vertex in $D_1(G) - V(P)$, it follows that $n_2 \le n  - N_3 - n_1 \le n-4$.
We conclude that $n-2t\leq\alpha_{k-reg}(T_n)\leq n-4$.

\noindent
(ii)
Suppose $\frac{n}{3}\leq t\leq n-5$. We first assume $n_1\geq \frac{n+2}{3}$.
Since $D_1(T_{n,t})$ is a $k$-independent set,
it follows that $\alpha_{k-reg}(T_{n,t})\geq n_1 \geq \frac{n+2}{3}$.
Next, we assume $n_1< \frac{n+2}{3}$, and so it follows
that $n_2\geq n-2n_1+2$, $n_1\geq N_3+2$ and $n_2=n-n_1-N_3$.
As $D_2(G)$ is a $k$-independent set, it follows
that $\alpha_{k-reg}(T_{n,t}) \ge  n_2 \ge  n-2n_1+2 \ge n-\frac{2n+4}{3}+2=\frac{n+2}{3}$.
The upper bound follows from (\ref{u-bound}).

\noindent (iii) Suppose that $t=n-4$. Then $diam(T_{n,t})=n-t=4$ and $h = 2$.
By (\ref{ff}), $f(n,t) =
\left\lceil\frac{n-1}{2}\right\rceil$.
Since vertices in $D_2(T)$ cannot be the end vertices of $P$ and of
$P_w$, for each $w \in D_1(T) - \{v_1, v_{n-t+1}\}$, and cannot be in $N_3$, it follows
that
\begin{eqnarray} \label{n2}
n_2 & \le & |V(P) - \{v_1, v_{j_0}, v_{n-t+1}\}| + \sum_{w \in D_1(T) - \{v_1, v_{n-t+1}\}}
(|V(P_w)| - 2)
\\ \nonumber
& = & (n-t+1) - 3 + (n_1 - 2)(h-1)
\end{eqnarray}
Since $n-t=4$ and $h = 2$, (\ref{n2}) leads to $n_2 \le n_1$.
By Lemma \ref{max=}, $\alpha_{k-reg}(T_{n,t}) = n_1$.
By Lemma \ref{n1-bound}, $t+1 \ge n_1 \ge f(n,t)= \left\lceil\frac{n-1}{2}\right\rceil$.
Thus (iii) must hold.

\noindent (iv) Suppose that $t=n-3$. Then $diam(T_{n,t})=n-t=3$ and
$h = 1$. Thus by (\ref{n2}), $n_2 < n_1$ and so by Lemma \ref{max=},
$\alpha_{k-reg}(T_{n,t}) = n_1$. By Lemma \ref{n1-bound}, $t+1 \ge
n_1 \ge f(n,t)$. By (\ref{ff}) with $n-t=3$, we have $f(n,t) = t+1$.
This implies (iv).

\qed
\end{pf}

The bounds in Theorem \ref{th3-1} are best possible in some sense,
as can be seen in the following examples.

\noindent {\bf Example 1:} $(1)$ Let $P = v_1\cdots v_{n-t+1}$ be a
path. For the lower bound,  let $v_1', v_2', \cdots, v_{t-1}'$ be
vertices not in $V(P)$ with $2 \le t \le \frac{n-1}{3}$. Since $n
\ge 2t$, there exists distinct vertices $v_{i_1}, v_{i_2}, \cdots,
v_{i_{t-1}} \in V(P) - \{v_1, v_{n-t+1}\}$. Obtain a $T_{n,t}$ with
$V(T_{n,t}) = V(P) \cup \{v_1', v_2', \cdots, v_{t-1}'\}$ and
$E(T_{n,t}) =  E(P) \cup \{v_{i_j}v_{j}': 1 \le j \le t-1\}$, (see
Figure 1 $(a)$). Then in this $T_{n,t}$, we have $n_1=t+1$,
$n_2=n-2t$, $N_3=t-1$. Thus any $k$-regular independent set $W$ must
be a subset of $D_j(G)$, for some $j$ with $1 \le j \le 3$. Since $n
\ge 3t+1$, we have $n_2 \ge n_1$. As $k \ge 2$, $D_2(G)$  is a
maximum regular $k$-independent set of $T_{n,t}$, and so
$\alpha_{k-reg}(T_{n,t}) = n - 2t$ if $n \ge 3t+1$.

For the upper bound, let $L = v_0'v_1'v_2'\cdots v_{t-1}'$ denote a
path. Obtain a $T_{n,t}'$ from $P$ and $T$ by identifying the vertex
$v_{j} \in V(P)$ and $v_0' \in V(L)$, where $j = \lfloor
\frac{n-t}{2} \rfloor+1$, (see Figure 1 $(b)$). In this case, we
have $n_1=3$, $n_2=n-4$, $N_3=1$, and so when $k \ge 2$ and $n \ge 7$,
$\alpha_{k-reg}(T_{n,t}')=n-4$, which shows the upper bound is sharp.

\begin{figure}[!hbpt]
\begin{center}
\includegraphics[scale=0.75]{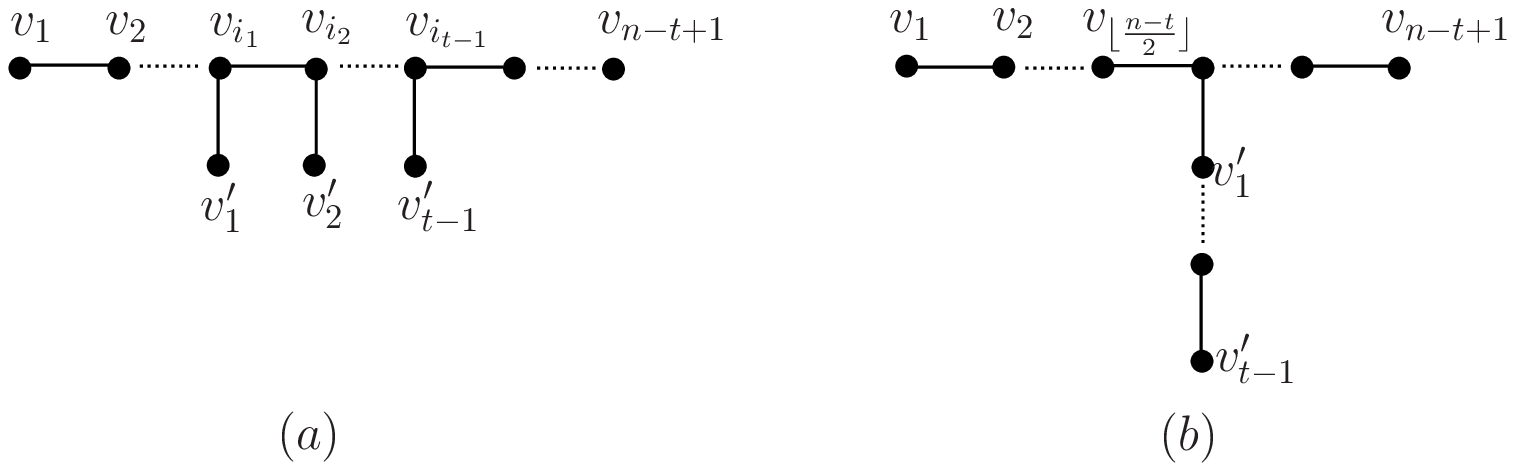}\\
Figure $1$. $(a)$ Tree with $\alpha_{k-reg}(T)=n-2t$ for $k\geq2$.
$(b)$ Tree with $\alpha_{k-reg}(T)=n-4$ for $k\geq2$.
\end{center}\label{fig1}
\end{figure}

\noindent {\bf Example 2:} $(2)$ For the lower bound, we let
$n=3t-2$ for some integer $t\ge 4$. Then $n-t+1 = 2t-1$. Let $P =
v_1v_2\cdots v_{2t-1}$ be a path. let $v_1', v_2', \cdots, v_{t-1}'$
be vertices not in $V(P)$. Obtain a $T_{n,t}(1)$ with $V(T_{n,t}(1))
= V(P) \cup \{v_1', v_2', \cdots, v_{t-1}'\}$ and  $E(T_{n,t}(1)) =
E(P) \cup \{v_{j+1}v_{j}': 1 \le j \le t-2\} \cup
\{v_{t-2}'v_{t-1}'\}$, (see Figure 2 $(a)$). Then in this
$T_{n,t}(1)$, we have $n_1=t=n_2$, $N_3=t-2$. As $k \ge 2$, each of
$D_1(T_{n,t}(1))$ and $D_2(T_{n,t}(1))$ is a maximum $k$-regular
independent set, and so
$\alpha_{k-reg}(T_{n,t}(1))=t=\frac{n+2}{3}$.

For the upper bound, let $n=12$ and $t=4$,
and let $T_{n,t}(2)$ be the tree depicted in
Figure 2(b). Then we have $n_1=3$, $n_2=8$, $N_3=1$.
It is routine to see that
$\alpha_{k-reg}(T_{n,t}(2))=n-f(n,t)-1=8$ for $k\geq2$.
Let $n=11$ and $t=4$,
and let $T_{n,t}(3)$ be the tree depicted in
Figure 2(c). Then we have $n_1=3$, $n_2=7$, $N_3=1$ and
$\alpha_{k-reg}(T_{n,t}(3))=n-f(n,t)-1=7$ for $k\geq2$.
Let $n=9$ and $t=4$, and let $T_{n,t}(4)$ be the tree depicted in
Figure 2(d). Then we have $n_1=5$, $n_2=1$, $N_3=3$ and
$\alpha_{k-reg}(T_{n,t}(4))=t+1=5$ for $k\geq2$.

\begin{figure}[!hbpt]
\begin{center}
\includegraphics[scale=0.75]{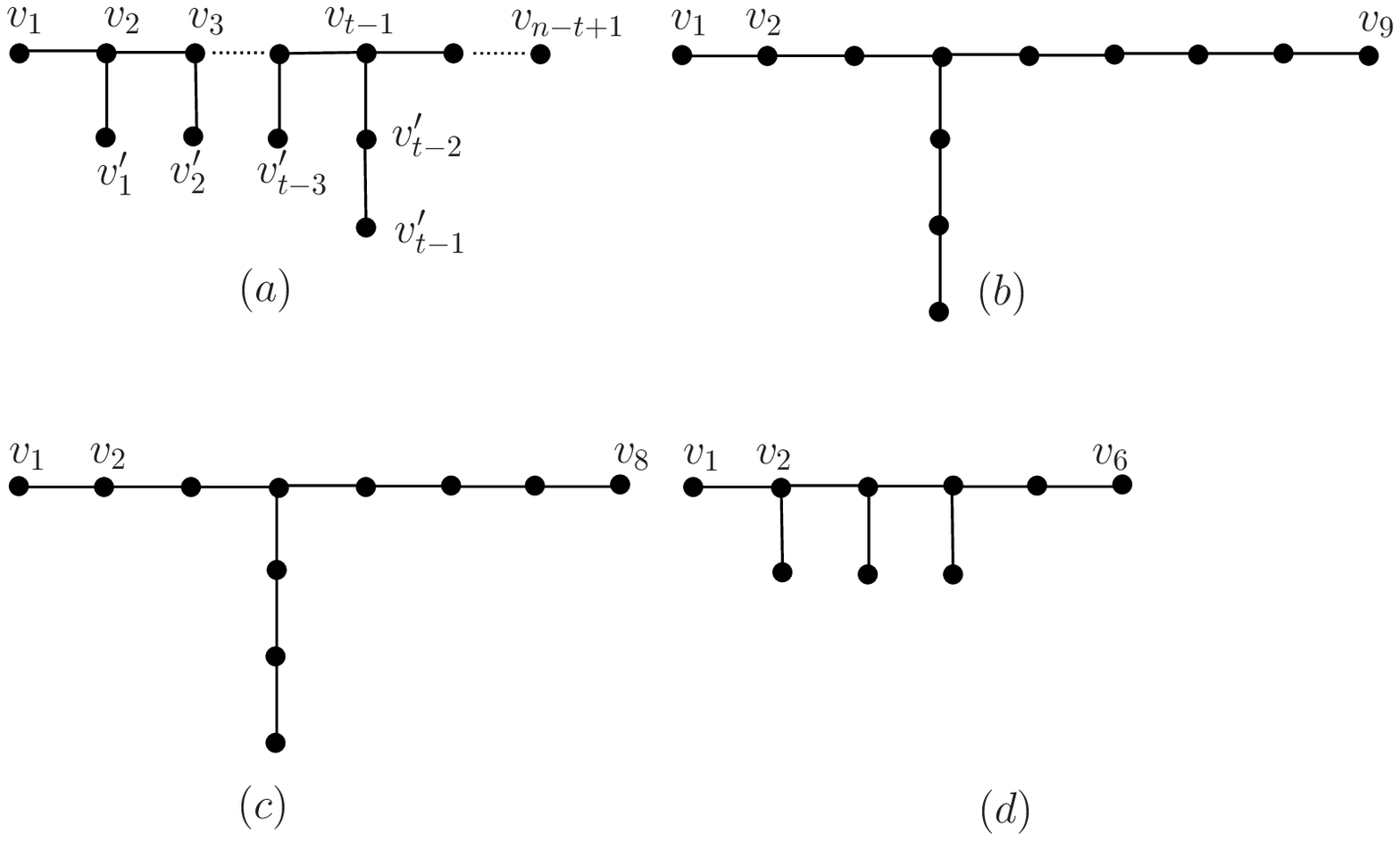}\\
Figure $2$. $(a)$ Tree with $\alpha_{k-reg}(T)=\frac{n+2}{3}$ for $k\geq2$.\\$(b)$ Tree with $\alpha_{k-reg}(T)=n-f(n,t)-1=8$ for $k\geq2$.\\ $(c)$ Tree with $\alpha_{k-reg}(T)=n-f(n,t)-1=7$ for $k\geq2$.\\
$(d)$ Tree with $\alpha_{k-reg}(T)=t+1=5$ for $k\geq2$.
\end{center}\label{fig2}
\end{figure}

\noindent {\bf Example 3:} $(3)$ For the lower bound,
let $n=10$ and $t=6$, ande let
$T_{n,t}(5)$ be the tre depicted in
Figure 3(a). Then we have
$n_1=5$, $n_2=4$, $N_3=1$,
and $\alpha_{k-reg}(T_{n,t}(5))=\left\lceil\frac{n-1}{2}\right\rceil=5$ for $k\geq2$.
For the upper bound, let $n=10$
and $t=6$, and  let
$T_{n,t}(6)$ be the tree depicted in
Figure 3(b). Then we have  $n_1=7$,
$n_2=2$, $N_3=1$ and $\alpha_{k-reg}(T_{n,t}(6))=t+1=7$ for $k\geq2$.

\begin{figure}[!hbpt]
\begin{center}
\includegraphics[scale=0.75]{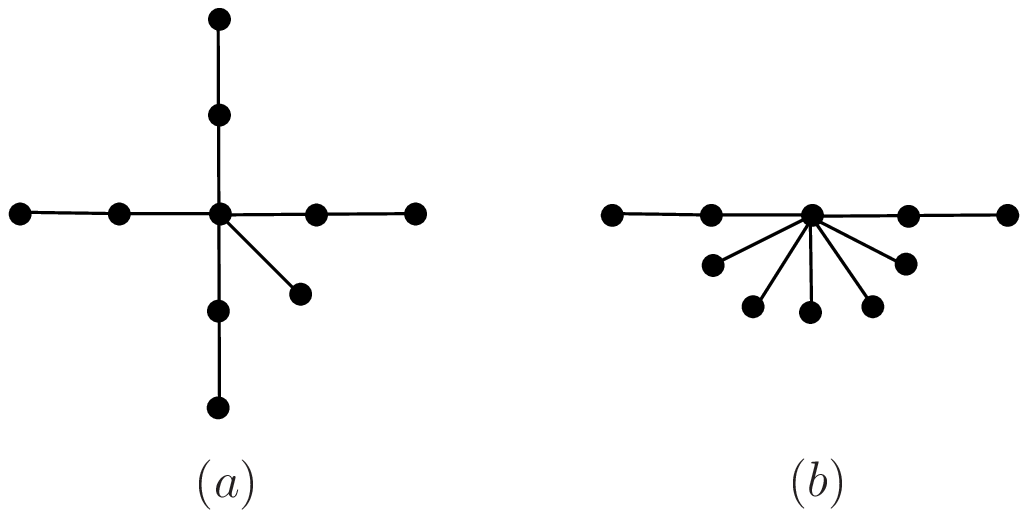}\\
Figure $3$.$(a)$ Tree with $\alpha_{k-reg}(T)=\left\lceil\frac{n-1}{2}\right\rceil=5$ for $k\geq2$.\\
$(b)$ Tree with $\alpha_{k-reg}(T)=t+1=7$ for $k\geq2$.
\end{center}\label{fig3}
\end{figure}

\section{Results for line graphs}

In this section, we investigate the bounds for the regular
$k$-independence number of line graphs
of graphs in certain families, including
trees, maximal outerplanar graphs and triangulations.
Throughout this section, $G$ denotes a graph with $m = |E(G)|$, and
define
$\delta' = \delta(L(G))$ and $\Delta' = \Delta(L(G))$.
For each integer $i$ with $\delta'\leq i\leq \Delta'$, define
$\Gamma_i = L(G)[D_i(L(G))]$.
%Use \Gamma_i for L(G_i) in the original paper.
Recall that the  repetition number of a graph $G$, defined in
(\ref{rep}), is the maximum $|D_i(G)|$ over all possible values of $i$.

\begin{lem} (Caro and Wes {\upshape \cite{Wes}})  \label{lem5-1}
Let $G$ be a graph with $m$ edges. Then
$$
rep(L(G))\geq \frac{1}{4}m^{\frac{1}{3}}.
$$
\end{lem}

\begin{thm}\label{th5-1}
Let $G$ be a graph with $m$ edges. Then
$$
\alpha_{k-reg}(L(G))\geq \frac{m^{\frac{1}{3}}}{4\chi_k(L(G))}.
$$
\end{thm}

\begin{pf}
Since $\alpha_{k,j}(L(G))\geq \frac{V(\Gamma_j)}{\chi_k(\Gamma_j)}$
holds for every $j$,
we have
\begin{equation} \label{a-lg}
 \alpha_{k-reg}(L(G))=\max \left\{\alpha_{k,j}(L(G)), \delta'\leq j \leq\Delta'\right \}\geq
 \max \left \{\frac{V(\Gamma_j)}{\chi_k(\Gamma_j)}: \delta'\leq j\leq \Delta' \right \}.
\end{equation}

Since $\chi_k(L(G))\geq \chi_k(L(G_j))$ holds for every $j$, it follows that
\begin{eqnarray} \label{ratio}
\max_{\delta'\leq j\leq \Delta'} \left\{\frac{V(\Gamma_j)}{\chi_k(\Gamma_j)} \right\}
& \geq &
 \max_{\delta'\leq j\leq \Delta'}
 \left\{\frac{rep(L(G))}{\chi_k(\Gamma_j)}: |V(\Gamma_j)|=rep(L(G))  \delta'\leq j\leq \Delta' \right\}
\\ \nonumber
& \geq & \frac{rep(L(G))}{\chi_k(L(G))}.
\end{eqnarray}
By Lemma \ref{lem5-1}, we have $rep(L(G))\geq \frac{1}{4}m^{\frac{1}{3}}$, and so
$\alpha_{k-reg}(L(G))\geq \frac{m^{\frac{1}{3}}}{4\chi_k(L(G))}$.

\qed
\end{pf}

\begin{lem} (Caro and Wes {\upshape \cite{Wes}}) \label{lem5-2}
Let $G$ be a graph with average degree $d$, minimum degree $\delta$, and $m$ edges.
If $d\geq \delta\geq1$, then
$$
rep(L(G))\geq \alpha\sqrt{m}-1,
$$
where $\alpha=\frac{\delta}{\sqrt{cd(cd-\delta)}}$ with $c=2d-2\delta+1$.
\end{lem}

\begin{thm}\label{th5-2}
Let $G$ be a graph with average degree $d$, minimum degree $\delta$, and $m$ edges. If $d\geq \delta\geq1$, then
$$
\alpha_{k-reg}(L(G))\geq \frac{\alpha\sqrt{m}-1}{\chi_k(L(G))},
$$
where $\alpha=\frac{\delta}{\sqrt{cd(cd-\delta)}}$ with $c=2d-2\delta+1$.
\end{thm}

\begin{pf}

Since for every $j$, we have $\alpha_{k,j}(L(G))\geq \frac{V(\Gamma_j)}{\chi_k(\Gamma_j)}$,
it follows that (\ref{a-lg}) must hold.
Hence, since for every $j$, we have $\chi_k(L(G))\geq \chi_k(\Gamma_j)$,
(\ref{ratio}) also holds.
By Lemma \ref{lem5-2}, $rep(L(G))\geq \alpha\sqrt{m}-1$,
where $\alpha=\frac{\delta}{\sqrt{cd(cd-\delta)}}$ with $c=2d-2\delta+1$.
It follows that $\alpha_{k-reg}(L(G))\geq \frac{\alpha\sqrt{m}-1}{\chi_k(L(G))}$.
\qed
\end{pf}

\begin{lem} (Caro and Wes {\upshape \cite{Wes}}) \label{lem5-3}
For a tree or maximal planar graph with $m$ edges, the repetition number of the line graph is at least $\sqrt{\frac{m}{30}}$ or $\sqrt{\frac{m}{182}}$, respectively.
\end{lem}

\begin{cor}\label{cor5-1}
$(1)$ If $G$ is a tree, then
$$
\alpha_{k-reg}(L(G))\geq \frac{\sqrt{m}}{\sqrt{30}\chi_k(L(G))}.
$$

$(2)$ If $G$ is a maximal planar graph with $m$ edges, then
$$
\alpha_{k-reg}(L(G))\geq \frac{\sqrt{m}}{\sqrt{182}\chi_k(L(G))}.
$$
\end{cor}

\begin{pf}
Let $n = |V(G)|$. If $G$ is a tree, then
$\delta=1$, $d=2-\frac{2}{n}$, and $c=3-\frac{4}{n}$,
It follows from Theorem \ref{th5-2} that
$\alpha_{k-reg}(L(G))\geq \frac{\sqrt{m}}{\sqrt{30}\chi_k(L(G))}$.
If $G$ is a maximal planar graph,
then $\delta=3$, $d=6-\frac{12}{n}$, and $c=7-\frac{12}{n}$.
It follows from Theorem \ref{th5-2} that
$\alpha_{k-reg}(L(G))\geq \frac{\sqrt{m}}{\sqrt{182}\chi_k(L(G))}$.
\qed
\end{pf}

\begin{lem}  (Caro and Wes {\upshape \cite{Wes}}) \label{lem5-4}
Let $G$ be a graph with $m$ edges. If $G$ is a tree with
a
perfect matching, a maximal outerplanar graph, or a triangulation with a $2$-factor, then $rep(L(G))$ is at
least $\frac{m}{6}$, $\frac{m}{14}$, or $\frac{m}{33}$, respectively. The lower bound improves to $\frac{m}{27}$ or $\frac{m}{15}$ for triangulations having a $2$-factor and minimum degree $4$ or $5$, respectively.
\end{lem}

\begin{thm}\label{th5-3}

$(1)$ If $G$ is a tree with a perfect matching, then
$$
\alpha_{k-reg}(L(G))\geq \frac{m}{6\chi_k(L(G))}.
$$

$(2)$ If $G$ is a maximal outerplanar graph with a $2$-factor, then
$$
\alpha_{k-reg}(L(G))\geq \frac{m}{14\chi_k(L(G))}.
$$

$(3)$ If $G$ is a triangulation graph with a $2$-factor, then
$$
\alpha_{k-reg}(L(G))\geq \frac{m}{33\chi_k(L(G))}.
$$

$(4)$ Moreover, if $G$ is a triangulation graph with a $2$-factor and minimum degree $4$ , then
$$
\alpha_{k-reg}(L(G))\geq \frac{m}{27\chi_k(L(G))}.
$$

$(5)$ Moreover, if $G$ is a triangulation graph with a $2$-factor and minimum degree $5$ , then
$$
\alpha_{k-reg}(L(G))\geq \frac{m}{15\chi_k(L(G))}.
$$
\end{thm}

\begin{pf}

Since for every $j$, we have $\alpha_{k,j}(L(G))\geq \frac{V(\Gamma_j)}{\chi_k(\Gamma_j)}$,
it follows that (\ref{a-lg}) must hold.
Hence, since for every $j$, we have $\chi_k(L(G))\geq \chi_k(\Gamma_j)$,
(\ref{ratio}) also holds.
By Lemma \ref{lem5-4},
\[
rep(L(G)) \ge
\left\{
\begin{array}{ll}
\frac{m}{6\chi_k(L(G))} & \mbox{ if $G$ is a tree with a perfect matching}
\\
\frac{m}{14\chi_k(L(G))} & \mbox{ if $G$ is  a maximal outerplanar graph}
\\
\frac{m}{33\chi_k(L(G))}  & \mbox{ if $G$ is  a triangulation with a $2$-factor}
\\
\frac{m}{27\chi_k(L(G))}  & \mbox{ if $G$ is  a triangulation with a $2$-factor with $\delta(G) \ge 4$}
\\
\frac{m}{15\chi_k(L(G))}  & \mbox{ if $G$ is  a triangulation with a $2$-factor with $\delta(G) \ge 5$}
\end{array}
\right. .
\]
Thus the conclusions of the theorem follows from (\ref{a-lg}) and (\ref{ratio}).

\qed
\end{pf}

\begin{lem} (Caro and Wes {\upshape \cite{Wes}}) \label{lem5-5}
Let $G$ be a triangulation with $m$ edges. If $G$ has minimum degree at least $4$, then $rep(L(G))\geq \frac{m}{68\chi_k(L(G))}$. If $G$ has minimum degree at least $5$, then $rep(L(G))\geq \frac{m}{51\chi_k(L(G))}$.
\end{lem}

\begin{thm}\label{th5-4}

$(1)$ If $G$ is a triangulation graph with $m$ edges and minimum degree at least $4$, then
$$
\alpha_{k-reg}(L(G))\geq \frac{m}{68\chi_k(L(G))}.
$$

$(2)$ If $G$ is a triangulation graph with $m$ edges and minimum degree at least $5$, then
$$
\alpha_{k-reg}(L(G))\geq \frac{m}{51\chi_k(L(G))}.
$$
\end{thm}

\begin{pf}

The proof of this theorem is similar to that of Theorem \ref{th5-3},
using Lemma \ref{lem5-5} instead of Lemma \ref{lem5-4}.
Therefore, it is omitted.

\qed
\end{pf}

\section{Graphs with given regular $k$-independence number}

The following proposition follows immediately from definition.
\begin{pro}\label{pro4-1}
Let $G$ be a simple graph of order $n$. Then
$$
1\leq\alpha_{k-reg}(G)\leq n.
$$
\end{pro}
In the rest of this section, we will present
characterizations for graphs reaching either bounds in
Proposition \ref{pro4-1}.

\begin{lem}\label{lem4-1}
Let $G$ be a simple graph with $n \ge 2$ vertices.
Then there exist at least two vertices with the same degree in $G$.
\end{lem}
\begin{pf}
This follows from the observation that either $G$ has an isolated vertex, and so
for any $v \in V(G)$, $0 \le d_G(v) \le n-2$;
or $G$ has no isolated vertices, and so
for any $v \in V(G)$, $1 \le d_G(v) \le n-1$.
\qed
\end{pf}

\begin{thm}\label{th4-1}
Let $G$ be a simple graph. Then $\alpha_{k-reg}(G)=1$ if and only
if $k=0$ and any subset of vertices with same degree in $G$ induces a clique of $G$.
\end{thm}
\begin{pf}
Suppose $\alpha_{k-reg}(G)=1$. Clearly, $k=0$ and $\alpha_{0-reg}(G)=1$.
Then $\alpha_{0-reg}(G)=\alpha_{reg}(G)=1$. By Lemma \ref{lem4-1},
there exist at least two vertices with same degree in $G$.
By the definition of the regular independence number,
there exists an edge between any two vertices with same
degree in $G$. Hence, any subset of vertices with same
degree in $G$ induces a clique in $G$, as desired.

Conversely, suppose that $k=0$ and any subset of vertices
with same degree in $G$ induces a clique of $G$. Note that
any two vertices with same degree in $G$ are adjacent.
By the definition of the regular independence number,
we have $\alpha_{reg}(G)=1$ and $\alpha_{0-reg}(G)=\alpha_{reg}(G)$.
So, $\alpha_{0-reg}(G)=1$.\qed
\end{pf}

\begin{thm}\label{th4-2}
Let $h$ be a nonnegative integer. Then $\alpha_{k-reg}(G)=n$
if and only if $G$ is a $h$-regular graph with $n$ vertices and $k\geq h$.
\end{thm}
\begin{pf}
Suppose $\alpha_{k-reg}(G)=n$. By the definition of the
regular $k$-independence number, we have that all vertices
in $G$ have same degree in $G$. Hence, the graph $G$ is a
$h$-regular graph with $n$ vertices and $k\geq h$.

Conversely, if $G$ is a $h$-regular graph with $n$ vertices
and $k\geq h$, then all vertices form a regular $k$-independent
set. By the definition of the regular $k$-independence number,
we have $\alpha_{k-reg}(G)=n$ for $k\geq h$. \qed
\end{pf}

\section{Nordhaus-Gaddum-type results}

In this section, we investigate the Nordhaus-Gaddum-type
problem on the regular $k$-independence number of graphs and obtain
the sharp bounds for
both
$\alpha_{k-reg}(G)+\alpha_{k-reg}(\Bar{G})$,
and
$\alpha_{k-reg}(G).\alpha_{k-reg}(\Bar{G})$,
over the class $\mathcal {G}(n)$ and characterize the extremal graphs.

\begin{thm}\label{th6-1}
For any $G\in \mathcal {G}(n)$, $\Bar{G}$ denotes the complement of $G$. Then

$(1)$ $3\leq \alpha_{k-reg}(G)+\alpha_{k-reg}(\Bar{G})\leq 2n$;

$(2)$ $2\leq\alpha_{k-reg}(G)\cdot \alpha_{k-reg}(\Bar{G})\leq n^2$.
\end{thm}
\begin{pf}
$(1)$ By Proposition \ref{pro4-1}, $\alpha_{k-reg}(G)\geq 1$
and $\alpha_{k-reg}(\bar{G})\geq 1$, and so
$\alpha_{k-reg}(G)+\alpha_{k-reg}(\Bar{G})\geq 2$.
However,
By Proposition \ref{pro4-1},
$\alpha_{k-reg}(G)+\alpha_{k-reg}(\Bar{G})=2$
if and only if $\alpha_{k-reg}(G)=1$ and $\alpha_{k-reg}(\bar{G})=1$.
Thus if $\alpha_{k-reg}(G)=1=\alpha_{k-reg}(\bar{G})$, then
by Theorem \ref{th4-1}, for any $i$,
either $D_i(G) = \emptyset$ or $G[D_i(G)]$ is a clique.
By the same reason,
either $D_{n-i-1}(\Bar{G}) = \emptyset$ or
$\Bar{G}[D_{n-i-1}(\Bar{G})]$
is a clique. Since $G \cup \Bar{G} = K_n$,
it is impossible that both $G[D_i(G)]$
and $\Bar{G})[D_{n-i-1}(\Bar{G})]$ are cliques.
This contradiction shows that we must have
$\alpha_{k-reg}(G)+\alpha_{k-reg}(\Bar{G})\geq 3$.
By Proposition \ref{pro4-1}, $\alpha_{k-reg}(G)\leq n$ and $\alpha_{k-reg}(\bar{G})\leq n$, and so $\alpha_{k-reg}(G)+ \alpha_{k-reg}(\Bar{G})\leq 2n$.

$(2)$
By Proposition \ref{pro4-1}, $\alpha_{k-reg}(G)\geq 1$ and $\alpha_{k-reg}(\bar{G})\geq 1$, and so $\alpha_{k-reg}(G)\cdot \alpha_{k-reg}(\Bar{G})\geq 1$.
By Proposition \ref{pro4-1},
$\alpha_{k-reg}(G)\cdot \alpha_{k-reg}(\Bar{G})=1$
if and only if $\alpha_{k-reg}(G)=1=\alpha_{k-reg}(\bar{G})$, and so
we obtain a contradiction as above.
Therefore, $\alpha_{k-reg}(G)\cdot \alpha_{k-reg}(\Bar{G})\geq 2$.
From Proposition \ref{pro4-1}, $\alpha_{k-reg}(G)\leq n$ and $\alpha_{k-reg}(\bar{G})\leq n$, and so $\alpha_{k-reg}(G)\cdot \alpha_{k-reg}(\Bar{G})\leq n^2$.

\qed
\end{pf}

 Before we study the graphs reaching the bounds
in Theorem \ref{th6-1}, we make the following observations.
\begin{obs} \label{obs1}
Let $n \ge 2$ be an integer, and $G\in \mathcal {G}(n)$. Then the following are equivalent.
\\
(i) $\alpha_{k-reg}(G)+\alpha_{k-reg}(\Bar{G})=3$.
\\
(ii) $\alpha_{k-reg}(G)\cdot \alpha_{k-reg}(\Bar{G})=2$.
\\
(iii) $\{\alpha_{k-reg}(G), \alpha_{k-reg}(\Bar{G})\} = \{1, 2\}$.
\end{obs}
In fact, by Proposition \ref{pro4-1},  each of (i) and (ii)
of Observation \ref{obs1} is equivalent to (iii).

\begin{obs} \label{obs2}
Let $n \ge 2$ be an integer, and  $G\in \mathcal {G}(n)$. Then the following are equivalent.
\\
(i) $\alpha_{k-reg}(G)+\alpha_{k-reg}(\Bar{G})=2n$.
\\
(ii) $\alpha_{k-reg}(G)\cdot \alpha_{k-reg}(\Bar{G})=n^2$.
\\
(iii) $\alpha_{k-reg}(G) = n = \alpha_{k-reg}(\Bar{G})$.
\end{obs}

\begin{pro}\label{pro6-1}
Let $n \ge 2$ be an integer, and let $G\in \mathcal {G}(n)$. Then
\begin{equation} \label{23}
\alpha_{k-reg}(G)+\alpha_{k-reg}(\Bar{G})=3
\mbox{ or } \alpha_{k-reg}(G)\cdot \alpha_{k-reg}(\Bar{G})=2
\end{equation}
if and only if $G$ satisfies the following conditions.
\\
$(1)$ $k=0$;
\\
$(2)$
for any integer $i \ge 0$, $|D_i(G)| \le 2$ and $|D_i(\bar{G})| \le 2$.
\\
$(3)$ Either every $G[D_i(G)]$ is connected for any $i$ such that
$D_i(G) \neq \emptyset$, or every $\Bar{G}[D_j(\Bar{G})]$ is
connected for any $j$ such that $D_j(\bar{G}) \neq \emptyset$.
\end{pro}
\begin{pf}
If (\ref{23}) holds, then by Observation \ref{obs1} and by symmetry,
we assume that $\alpha_{k-reg}(G)=1$ and
$\alpha_{k-reg}(\Bar{G})=2$. By Theorem \ref{th4-1}, we have $k=0$
(and so (1) holds), and for each $i$ with $|D_i(G)| \ge 2$,
$G[D_i(G)]$ must be a clique. If $|D_i(G)| \ge 3$, then $D_i(G)$ is
a regular independent set in $\Bar{G}$, contrary to the fact that
$\alpha_{k-reg}(\Bar{G})=2$. This implies that
\begin{equation} \label{Di}
|D_i(G)| \le 2, \mbox{ for every $i$.}
\end{equation}
If for some $j$, $|D_{j}(\Bar{G})| \ge 3$, then
as $D_{n-j-1}(G) = D_j(\Bar{G})$, we have $|D_{n-j-1}(G)| \ge 3$,
contrary to (\ref{Di}). Hence $|D_{j}(\Bar{G})| \le 2$ for
every $j$, and so (2) holds.

To justify (3), we observe that if for some $i \neq j$,
both $G[D_i(G)]$ and $\Bar{G}[D_j(\Bar{G})]$ are disconnected,
then by (2), $G[D_i(G)]$ is independent in $G$ and
$\Bar{G}[D_j(\Bar{G})]$ is independent in $\Bar{G}$.
It follows that $\alpha_{k-reg}(G)+\alpha_{k-reg}(\Bar{G})
\ge |D_i(G)| + |D_j(\Bar{G})| \ge 2 + 2 = 4$,
contrary to (\ref{23}). Thus (3) must hold.

Conversely, suppose $G$ satisfies Proposition \ref{pro6-1} (1), (2)
and (3). Then by $n \ge 2$ and by (2), every regular 0-independent
set in $G$ and in $\Bar{G}$ is of size at most 2. By (3),
$\{\alpha_{k-reg}(G), \alpha_{k-reg}(\Bar{G})\}=\{1,2\}$, and so by
Observation \ref{obs1}, (\ref{23}) must hold. \qed
\end{pf}

\noindent {\bf Remark 1:} In fact,
graphs satisfying the conditions
in Proposition \ref{pro6-1}
indeed exist. Let $H$ be a graph obtained from
a triangle and an edge by identifying a vertex of this triangle and an
endpoint of this edge. Clearly, $|V(H)|=4$. Let $G$ be a union of $H$
and an isolated vertex. Then, $|V(G)|=5$.
It is routine to check that
graphs $G$ and $\Bar{G}$ satisfy $\alpha_{k-reg}(G)+\alpha_{k-reg}(\Bar{G})=3$
and $\alpha_{k-reg}(G)\cdot \alpha_{k-reg}(\Bar{G})=2$ for $k=0$.

\begin{pro}\label{pro6-2}
Let $n \ge 2$ and $h \ge 0$ be integers, and
let $G\in \mathcal {G}(n)$.
The following are equivalent.
\\
(i) $\alpha_{k-reg}(G)+\alpha_{k-reg}(\Bar{G})=2n$.
\\
(ii) $\alpha_{k-reg}(G)\cdot \alpha_{k-reg}(\Bar{G})=n^2$.
\\
(iii) $G$ is an $h$-regular graph with $n$
vertices and $k\geq\max\{h,n-1-h\}$.
\end{pro}

\begin{pf}
By Observation \ref{obs2}, it suffices to show that (i) and (iii)
are equivalent.
%{\color{blue}it suffices to show that (i) and (ii) } are equivalent.

Assume that (i) holds. By Observation \ref{obs2},
$\alpha_{k-reg}(G)=n$, and so $V(G)$ is a regular $k$-independent set.
Hence  by Theorem \ref{th4-2}, $G$ must be $h$-regular graph and $k \ge h$.
By  Observation \ref{obs2},
$\alpha_{k-reg}(\Bar{G})=n$, and so the same argument shows that
$\Bar{G}$ is $n-1-h$-regular and $k\geq n-1-h$, and so (iii) must hold.

Assume that (iii) holds. Then $G$ is an $h$-regular graph of order $n$,
and $\Bar{G}$ is  $(n-1-h)$-regular.
Since $k\geq\max\{h,n-1-h\}$, $V(G)$ is a regular $k$-independent set
in both $G$ and $\Bar{G}$, and so  by Theorem \ref{th4-2},
$\alpha_{k-reg}(G)=n=\alpha_{k-reg}(\Bar{G})=n$. Thus (i) follows.
\qed
\end{pf}

\end{document}